\documentclass{amsart}

\usepackage{amsmath,amsthm,amsfonts,amssymb,bm,graphicx, mathrsfs}

\usepackage
{hyperref}
\hypersetup{colorlinks=true,citecolor=blue,linkcolor=blue,urlcolor=blue,
pdfstartview=FitH }

\theoremstyle{plain}
\newtheorem{theorem}{Theorem}

\newtheorem{lemma}{Lemma}

\newtheorem{cor}[theorem]{Corollary}
\newtheorem{prop}[theorem]{Proposition}
\numberwithin{equation}{section}

\theoremstyle{definition}

\renewcommand{\geq}{\geqslant}
\renewcommand{\leq}{\leqslant}

\newcommand{\changed}[1]{{\color{black} #1}}

\usepackage{calc}
\newsavebox\CBox
\newcommand\hcancel[2][0.5pt]{%
  \changed{\ifmmode\sbox\CBox{$#2$}\else\sbox\CBox{#2}\fi%
  \makebox[0pt][l]{\usebox\CBox}%
  \rule[0.5\ht\CBox-#1/2]{\wd\CBox}{#1}}}

\makeatletter
\DeclareRobustCommand\widecheck[1]{{\mathpalette\@widecheck{#1}}}
\def\@widecheck#1#2{%
    \setbox\z@\hbox{\m@th$#1#2$}%
    \setbox\tw@\hbox{\m@th$#1%
       \widehat{%
          \vrule\@width\z@\@height\ht\z@
          \vrule\@height\z@\@width\wd\z@}$}%
    \dp\tw@-\ht\z@
    \@tempdima\ht\z@ \advance\@tempdima2\ht\tw@ \divide\@tempdima\thr@@
    \setbox\tw@\hbox{%
       \raise\@tempdima\hbox{\scalebox{1}[-1]{\lower\@tempdima\box
\tw@}}}%
    {\ooalign{\box\tw@ \cr \box\z@}}}
\makeatother

\begin{document}

\author{Valentin Blomer}
  
\address{Mathematisches Institut, Bunsenstr. 3-5, 37073 G\"ottingen, Germany} \email{vblomer@math.uni-goettingen.de}

 \title{Epstein zeta-functions, subconvexity, and the purity conjecture}

\thanks{First author   supported in part  by the  Volkswagen Foundation and  NSF grant 1128155 while enjoying the hospitality of 
the Institute for Advanced Study. The United States Government is authorized to reproduce and distribute reprints notwithstanding any copyright notation herein.}


\begin{abstract}  Subconvexity bounds are proved for general Epstein zeta functions of $k$-ary quadratic forms. This is related to sup-norm bounds for Eisenstein series on ${\rm GL}(k)$, and the exact  sup-norm exponent is determined to be $(k-2)/8$ for $k \geq 4$. In particular, if $k$ is odd, this exponent is not in $\Bbb{Z}/4$, which shows that Sarnak's purity conjecture does not hold for Eisenstein series. 
 \end{abstract}


\setcounter{tocdepth}{2}  \maketitle 

\maketitle

\section{Introduction}

Eisenstein series are a good testing ground for properties of automorphic forms. In this note we are concerned with general Epstein zeta-functions. Let $k \geq 2$ be an integer, and let $Q(\textbf{x}) = \textbf{x}^{\top} Z \textbf{x}$ be a (not necessarily rational) positive-definite quadratic form in $k$ variables with symmetric matrix $Z$. The corresponding Epstein zeta-function \cite{Ep}
$$E(Z, s) := \sum_{\textbf{x} \in \Bbb{Z}^k \setminus \{0\}} Q(\textbf{x})^{-ks/2}$$
can be continued to the   complex plane with a simple pole at $s=1$ 
and satisfies the functional equation
\begin{equation}\label{fe}
  \Gamma_{\Bbb{R}}(ks) E(Z,  s) = (\det Z)^{-1/2}  \Gamma_{\Bbb{R}}(k(1-s)) E(Z^{-1}, 1-s)
  \end{equation}
with $\Gamma_{\Bbb{R}}(s) = \pi^{-s/2} \Gamma(s/2).$ Here and for the rest of the paper we normalize all Dirichlet series such that the ``critical strip'' is the region $0 \leq \Re s  \leq 1$. 
For a given $Z$, the function $s\mapsto E(Z,  s)$ is a generalized Dirichlet series that in most cases is not in the Selberg class, but nevertheless satisfies a functional equation essentially of ${\rm GL}(k)$-type. 

On the other hand, for fixed $s = 1/2 + it$,  the function $z \mapsto E(z^{\top} z, s)$ for $z \in  \Bbb{H}_k = {\rm GL}(k)/O(k) \Bbb{R}^{\ast}$ in the generalized upper half plane \cite{Go} is $\Gamma = {\rm SL}_k(\Bbb{Z})$-invariant and an eigenfunction of the ring of invariant differential operators of the locally symmetric space $\Gamma \backslash \Bbb{H}_k$. In fact, up to a constant factor $\zeta(k)$ it is a totally degenerate Eisenstein series associated to the maximal parabolic subgroup of ${\rm SL}_k(\Bbb{Z})$, see \cite[10.7]{Go}. It is therefore of interest to study its properties as an automorphic form. We will consider both viewpoints in this note.

\subsection{Subconvexity} The Dirichlet series $s \mapsto E(Z, s)$ is absolutely convergent in $\Re s > 1$, and the critical strip is $0 \leq \Re s \leq 1$. By the functional equation, the convexity bound on the critical line $\Re s = 1/2$ is 
$$E(Z, 1/2 + it) \ll_{\varepsilon, Z} (1+|t|)^{k/4+\varepsilon}.$$
 The subconvexity problem is one of the central topics for arithmetic Dirichlet series, but has so far only been solved for $L$-functions of degree 2 and in a handful of special cases of larger degree $\leq 8$. Our first result is a subconvexity bound for general Epstein $\zeta$-functions of arbitrary degree.
\begin{theorem}\label{thm1} Let $k \geq 2$. Let $\delta_2 = 1/6$, $\delta_3 = 1/4$  and $\delta_k = 1/2$ for $k \geq 4$. Then
\begin{equation}\label{sub}
E(Z, 1/2 + it) \ll_{\varepsilon, Z} (1+|t|)^{\frac{k}{4} - \delta_k+\varepsilon}
\end{equation}
for any $\varepsilon > 0$. The implied constant depends continuously on $Z$, in particular the bound is independent of $Z$ as long as $Z$ varies in a fixed compact domain. 
\end{theorem}
This improves in particular the value $\delta_2 = 1/8$, recently obtained by M.\ Young \cite{Yo} (which in turn improves the value $\delta_2 = 1/12$ of Iwaniec and Sarnak \cite{IS}). Numerically, going from $1/8$ to $1/6$ is the typical difference between a Burgess-type exponent coming from amplification and a stronger  Weyl-type bound coming from van der Corput-type estimates. For $k \geq 3$, no bounds of type \eqref{sub} seem to be in the literature. Fomenko \cite{Fo} has estimated $E(Z, s)$ for integral forms $Z$ on the line $\Re s = 1 -1/k$. It is not unlikely that Theorem \ref{thm1} remains true with $\delta_k = 1/2$ for all $k \geq 2$. 
The proof of Theorem \ref{thm1} relies on results of G\"otze \cite{Goe} for $k \geq 3$ and a method of Titchmarsh for $k=2$. \\

In general, Theorem \ref{thm1}  is (up to the value of $\varepsilon$) best possible for $k \geq 4$. 
\begin{theorem}\label{thm2} Let $k \geq 4$, and let $I_k$ be the $k$-by-$k$ identity matrix. Then
$$E(I_k, 1/2 + it) = \Omega\left( (1+|t|)^{\frac{k}{4} - \frac{1}{2}}\right).$$
\end{theorem}
In particular, Epstein-zeta functions of degree $k \geq 4$  in general do not satisfy the Lindel\"of hypothesis. (They may also violate the Riemann hypothesis quite badly, for instance $E(I_4, s) = 8(1 - 2^{1-2s}) \zeta(2s)\zeta(1 - 2s)$ has \emph{no} zeros on the critical line.) The proof of Theorem \ref{thm2}  starts with an   application of Siegel's mass formula. For odd $k$ -- and we will see in the next subsection that this is the most interesting case -- things are a little more complicated, and $E(I_k, 1/2 + it)$ turns out to be essentially  a special value of a double Dirichlet series outside the region of absolute convergence.  It requires   several applications of the various functional equations to bound this quantity from below.

The reason why the Eisenstein series at the identity exhibit such strong growth on the critical line is a rigidity phenomenon:  the values of integral forms are restricted to integers. In fact, with a little more work one can show the lower bound of Theorem \ref{thm2} for any fixed rational matrix. The situation changes dramatically for ``generic''   quadratic forms. It follows essentially from the Maa{\ss}-Selberg relations  to  compute the inner product of truncated Eisenstein series for ${\rm GL}(k)$ \cite[Lemma 4.2]{Ar} that for a fixed compact domain $\Omega \subseteq {\rm SL}_k(\Bbb{Z}) \backslash \Bbb{H}_k$ we have
\begin{equation}\label{maassselberg}
\int_T^{T+1} \int_{\Omega} |E(z ^{\top} z, 1/2 + it)|^2  d\mu(z)\, dt \ll_{\Omega, \varepsilon} T^{\varepsilon}
\end{equation}
(where $d\mu(z)$ is the usual hyperbolic measure on $\Bbb{H}_k$). This is worked out for $k = 3$ in complete detail in \cite{Mill} (in which case  it follows also from the stronger QUE result in \cite{Zh}).  From \eqref{maassselberg} one can easily conclude 
$$\int_0^{X} |E(Z, 1/2 + it)|^2 dt \ll X^{1+\varepsilon}, \quad E(Z, 1/2 + it) \ll (1+|t|)^{1/2 + \varepsilon}$$
for almost all $Z \in \Omega$. If $k=2$, Spinu \cite{Sp} has obtained bounds of the same quality as in \eqref{maassselberg} for the fourth moment, which implies $E(Z, 1/2 + it) \ll (1+|t|)^{1/4 + \varepsilon}$ for $k=2$ and almost all $Z \in \Omega$. 

We do not go into further details here, but consider  the variation of an ``almost all'' result within  
\emph{diagonal} forms. 
This uses a method of Jarn\'ik \cite{Ja}. 

\begin{prop}\label{thm3} Let $k \geq 2$. Let $Z$ vary in a fixed compact domain $D$ of positive diagonal matrices. Then there is a set $M \subseteq D$ of full measure, such that all for all $Z \in M$ we have
\begin{equation}\label{meansquare}
\int_0^{X} |E(Z, 1/2 + it)|^2 dt \ll X^{1+\varepsilon},
\end{equation}
in particular 
\begin{equation}\label{indiv}
E(Z, 1/2 + it) \ll (1+|t|)^{1/2 + \varepsilon}.
\end{equation}
\end{prop} 

In other words, if we define $C(t) := (1 + |t|)^{k}$ to be the analytic conductor of $E(Z, 1/2 + it)$, then we obtain the individual bound $E(Z, 1/2 + it) \ll C(t)^{1/(2k) + \varepsilon}$ for generic diagonal quadratic forms, which for large $k$ is almost as strong as the Lindel\"of hypothesis. The set $M$ is described explicitly in diophantine terms in \cite{Ja}. 

\subsection{Sup-norms of automorphic forms} In this subsection we re-interpret the subconvexity results in an automorphic context. The distribution of mass of an eigenfunction of the ring of differential operators on a symmetric space has received a lot of attention in the context of quantum chaos. It is a central question to what extent high energy eigenfunctions $\phi$, i.e.\ having a large Laplacian eigenvalue $\lambda$,  behave like random waves or display some structure related to the classical trajectories. One measure of equidistribution is a bound for the sup-norm $\|\phi \|_{\infty}$. While the random wave model predicts that $\| \phi \|_{\infty}$ cannot grow too quickly with $\lambda$ and certainly 
\begin{equation}\label{small}
   \| \phi \|_{\infty} \ll \lambda^{\varepsilon},
   \end{equation}
there are known phenomena where   \eqref{small} is violated. One is the behaviour close to the cusps \cite{Sa, BT} if the underlying space is not compact. This is an artefact of special functions  that we do not focus on in this discussion. Even in a compact part of the manifold, however, large sup-norms can occur when $\phi$ is a lift coming from a smaller group. This was first observed by Rudnick and Sarnak \cite{RS} for $\theta$-lifts on arithmetic 3-manifolds, based on similar phenomena of Eisenstein series on ${\rm SO}(3, 1)$. This has been extended and generalized in \cite{Do, Mili}; in particular, for each $n \geq 5$ there exist compact $n$-dimensional hyperbolic manifolds of constant negative curvature having a sequence of $L^2$-normalized eigenfunctions $\phi_j$ with 
\begin{equation}\label{4}
  \|\phi_j \|_{\infty} \gg \lambda_j^{(n-4)/4}.
 \end{equation}

In this context, Sarnak \cite{Sa} made the \emph{purity conjecture} that the set of accumulation points of $$\frac{\log \| \phi_j \|_{\infty}}{\log \lambda_j}$$
is contained in $\Bbb{Z}/4$. Up until now, this is consistent with all known examples and in particular with \eqref{4}, but it should be pointed out, however, that not a single non-zero  accumulation point on a negatively curved manifold has ever been determined explicitly. 

Theorems \ref{thm1} and \ref{thm2} can be interpreted in this context and  show that the obvious variation of Sarnak's purity  conjecture for Eisenstein series does \emph{not} hold. The (totally degenerate) Eisenstein series
$$ \frac{1}{\zeta(k)} E(z^{\top} z, s) = \sum_{\gamma \in P_{\text{max}}\backslash {\rm SL}_k(\Bbb{Z})} \det(\gamma z)^s$$
for $s = 1/2 + it$ (after analytic continuation) is an eigenfunction of the ring of invariant differential operators on ${\rm SL}_k(\Bbb{Z}) \backslash \Bbb{H}_k$ with Langlands parameters $(t, \ldots, t, - (k-1)t)$ and Laplacian eigenvalue  
$$\lambda = \lambda(t)  =  \frac{k^3 - k}{24} + \frac{1}{2}\bigl(t^2 + \ldots + t^2  + (k-1)^2 t^2\bigr) \asymp t^2.$$
It is, of course, not $L^2$-integrable on the whole space, but it is in $L^2$ on each   compact  subset. 
\begin{cor} Let $k \geq 4$ and let $\Omega \subseteq {\rm SL}_k(\Bbb{Z}) \backslash \Bbb{H}_k$ be compact.   Then the function
$$t \mapsto \frac{\log \| E(., 1/2 + it)|_{\Omega} \|_{\infty}}{\log \lambda(t)}$$
has the accumulation point $k/8 - 1/4$. In particular, if $k$ is odd, this is in $\Bbb{Z}/8$, but not in $\Bbb{Z}/4$. 
\end{cor}

\section{Proof of Theorem \ref{thm1}}\label{sec2} Let $|t| \geq 10$ and let $k \geq 2$ be fixed. We start with an approximate functional equation for $E(Z, s)$ on the critical line. A Hardy-Littlewood style version can be found in \cite[Theorem 2]{CN},  
but it is more convenient to use a smooth version. By a slight modification of   \cite[Theorem 5.3, Proposition 5.4]{IK} we conclude from the functional equation \eqref{fe} that 
\begin{equation}\label{IK}
E(Z, 1/2 + it) = \sum_{\textbf{x} \not= 0} \frac{W_t^+(Q_+(\textbf{x})^{k/2})}{Q_+(\textbf{x})^{k/2(1/2 + it)}} + \frac{\Gamma_{\Bbb{R}}(k(1/2 -  i t))}{\Gamma_{\Bbb{R}}(k(1/2  +  i t)) (\det Z)^{1/2}} \sum_{\textbf{x} \not= 0} \frac{W_t^-(Q_-(\textbf{x})^{k/2})}{Q_-(\textbf{x})^{k/2(1/2 - it)}} + O(1)
\end{equation}
where $Q_+ = Q$, $Q_-(\textbf{x}) = \textbf{x}^{\top} Z^{-1} \textbf{x}$ and 
$$W^{\pm}(y) =  \frac{1}{2\pi i} \int_{(1)} e^{u^2} 
\frac{\Gamma_{\Bbb{R}}(k(u + 1/2 \pm it))}{\Gamma_{\Bbb{R}}(k(1/2 \pm it) )}   y^{-u}  \frac{du}{u}.  $$
 The error term is a (crude) bound for the contribution of the pole of $E(Z, s)$ at $s=1$, and the exponential factor $e^{u^2}$ makes the integral rapidly convergent. Shifting the contour to the far right, we see that $W^{\pm}(y)$ becomes negligible, once $y \geq (1+|t|)^{k/2 + \varepsilon}$, so we can truncate the sums at $Q(\textbf{x})^{k/2} \leq |t|^{k/2 + \varepsilon}$ at the cost of a negligible error.  We attach a smooth partition  of unity to the $\textbf{x}$-sum that localizes at   $Q(\textbf{x}) \asymp T$, shift the contour back to $\Re u = \varepsilon$, truncate it  at $v = \Im u\ll |t|^{\varepsilon}$  and interchange summation and integration to obtain
\begin{equation}\label{afe}
E(Z, 1/2 + it)  \ll 1 + |t|^{\varepsilon} \sum_{\pm}   \sum_{ T = 2^{\nu}  \ll |t|^{1+\varepsilon}}  \int_{|v| \leq |t|^{\varepsilon} } T^{-k/4} \Biggl|\sum_{\textbf{x} \not= 0} \frac{V_{T}(Q_{\pm}(\textbf{x}))}{Q_{\pm}(\textbf{x})^{\pm ikt/2 + iv}} \Biggr| dv,
\end{equation}
where $T$ runs over powers of 2 and  $V_{T}$ has  compact support in $[T, 3T]$ and satisfies
\begin{equation}\label{v}
   \frac{\partial^j}{\partial x^j} V_{T}(x) \ll_{j} T^{-j}
   \end{equation}
for $j \in \Bbb{N}_0$.  Let
\begin{equation}\label{curlyV}
\mathcal{V}_{T, v}(s;t) := \int_0^{\infty} V_{T} (x) x^{\mp ikt/2 - iv} e^{s x} dx
\end{equation}
denote the Laplace transform (composed with $s \mapsto -s$). This is an entire function in $s$. Let
$$\theta(Q, s)  = \sum_{\textbf{x} \not= 0} e^{-s Q(\textbf{x})}.$$
By the Laplace inversion formula we have
\begin{equation}\label{laplace}
\sum_{\textbf{x} \not= 0} \frac{V_{ T}(Q_{\pm}(\textbf{x}))}{Q_{\pm}(\textbf{x})^{\pm ikt/2 + iv}}  
 =  \int_{(\delta)} \mathcal{V}_{T, v}(s;t) \theta(Q_{\pm}, s) \frac{ds}{2\pi i},
 \end{equation}
where $\delta > 0$ is arbitrary. In order to keep $ \mathcal{V}_{T, v}(s;t) $ manageable, we choose $\delta = 1/T$ and write $s = 1/T + i\tau$. It is not hard to evaluate $ \mathcal{V}_{T, v}(s;t) $ by a stationary phase argument, but we get an upper bound in a completely elementary fashion by the ``Gau{\ss} sum trick''. First we observe that by \eqref{v} and  partial integration $\mathcal{V}_{T, v}(s;t) \ll |t|^{-100}$ unless $\tau \asymp |t|/T$. This follows, for instance, from \cite[Lemma 8.1]{BKY} with $X = 1$, $U = Q = T$, $Y = |t|$ and $R = |\tau| + |t|/T$, provided $|\tau| \geq c |t|/T$ or $|t|/T \geq c |\tau|$ for a sufficiently large constant $c$.  
Under the condition $|\tau| \asymp |t|/T$, which we assume from now on, we have
\begin{displaymath}
\begin{split}
  |\mathcal{V}_{T, v}(s;t)|^2 & = \int_0^{\infty}\int_0^{\infty} V_T(x) \overline{V_T(y)}  \left(\frac{x}{y}\right)^{\pm ikt/2 + iv} e^{ (x+y)/T} e^{i\tau(x-y)} dx\, dy\\
  &  = \int_0^{\infty}\int_0^{\infty} V_T(xy) \overline{V_T(y)} x^{\pm ikt/2 + iv} e^{ (x+1)y/T} e^{i\tau(x-1)y} dx\, y\, dy.
\end{split}
\end{displaymath}
Again by \eqref{v} and partial integration, the $y$-integral is $\ll T^2 w( T\tau (x-1))$ for a rapidly decaying function $w$, so that by trivial estimates  
\begin{equation}\label{v1}
\mathcal{V}_{T, v}(s;t)  \ll \left(\frac{T^2}{T \tau}\right)^{1/2} \asymp \frac{T}{|t|^{1/2}}.
\end{equation} 

To bound the $\theta$-series in \eqref{laplace}, we use deep results of G\"otze.  By  \cite[Lemma 3.5, Lemma 3.10, (3.15)]{Goe} with
$$a, b \asymp \frac{|t|}{T}, \quad g(x) \ll 1, \quad G(a, b) \ll \frac{|t|}{T}, \quad d = k \geq 3, \quad r = T^{1/2}, \quad \gamma \geq 1, \quad D \ll T^{d/2}$$
(notice that for these results the assumption $d \geq 5$ in \cite{Goe} is not used and \cite[(3.15)]{Goe} implies $\gamma \geq  1$)  we have 
\begin{displaymath}
\begin{split}
  \int_{|\tau| \asymp |t|/T }   \left| \theta\left(Q_{\pm}, \frac{1}{T} + i\tau\right)\right| d\tau&  \ll T^{k/4} \Biggl(T^{k/4 - 1} \int_1^{T^{k/2}} v^{-1/2 + 2/k} \left(\frac{|t|}{T}  + 1\right)\frac{dv}{v} + \frac{|t|}{T}\Biggr)\\
  &\ll T^{k/4+\varepsilon} \frac{|t|}{T} \left( T^{k/4 - 1}(1 + \delta_{k=3} T^{1/4})+ 1 \right).
  \end{split}
  \end{displaymath} 
Substituting this together with \eqref{v1} and \eqref{laplace} into \eqref{afe}, we obtain
$E(Z, 1/2 + it) \ll |t|^{k/4 - 1/2 + \varepsilon}$ for $k \geq 4$ and $E(Z, 1/2 + it) \ll |t|^{1/2+\varepsilon}$ for $k = 3$. \\

For the remaining case $k =2$ we use a different strategy, which follows an argument of Titchmarsh \cite{Ti}, based on van der Corput's method. Titchmarsh showed 
$$E(Z, 1/2 + it) \ll |t|^{1/3} \log |t|$$
for $|t| \geq 2$ and \emph{diagonal} $Z$; here we make some minor simplifications of his argument and indicate briefly the necessary changes and additions for   arbitrary $Z$, given by a positive definite quadratic form $Q(\textbf{x}) = ax_1^2 + b x_1x_2 + cx_2^2$.  We may, without loss of generality, assume to be reduced, that is, $|b| \leq a \leq c$ (then also $Q_-$ is, up to permutation of $x_1$ and $x_2$, reduced), and in particular 
\begin{equation}\label{norm}
Q_{\pm}(\textbf{x}) \asymp \| \textbf{x} \|^2,
\end{equation}
 which we use frequently.  
On the one hand, we use a  more streamlined and effective approximate functional equation that makes \cite[Theorem 1]{Ti} and also the discussion of the boundary values at the end of the paper  obsolete. Secondly, we slightly simplify some of the estimations of multi-dimensional derivatives.  

With this in mind, we return to \eqref{IK}, specialize $k = 2$,  and remove the weight function $W_T^{\pm}(Q_{\pm}(\textbf{x}))$ as well as the factor $Q_{\pm}(\textbf{x})^{-1/2}$ by partial summation. Thus we need to prove
\begin{equation}\label{toprove}
\mathcal{S} := X^{-1} \sum_{ \substack{X_1 \leq x_1 \leq 2X_1\\ X_2 \leq x_2 \leq 2X_2}} e^{-i f(x_1, x_2)} 
\ll |t|^{1/3 + \varepsilon}, \quad f(x_1, x_2) = t\log Q_{\pm}(x_1, x_2), 
\end{equation}
for $X_1, X_2 \ll |t|^{1/2 + \varepsilon}$, $X := \max(X_1, X_2)$  and, without loss of generality,  $\min(X_1, X_2) \gg |t|^{1/3}$, for otherwise the bound is trivial. With the aim of applying Weyl differencing, we consider  exponential sums of the type
\begin{equation}\label{exp}
\mathcal{S}_1  = \mathcal{S}_1(\bm \mu):= \sum_{\substack{X_1 \leq x_1 \leq X_1' \\ X_2 \leq x_2 \leq X_2'}} e^{i g(x_1, x_2)}
\end{equation}
with $X_1' \leq X' \leq 2X_1$, $X_2' \leq X' \leq 2 X_2$ and 
$$g(x_1, x_2) =f(x_1, x_2) - f(x_1 + \mu_1, x_2 + \mu_2) =   t \log Q(x_1, y_1) - t \log Q(x_1 + \mu_1, x_2 + \mu_2)$$
with 
\begin{equation}\label{rho}
0 \not= \| {\bm \mu} \| \leq \rho := X |t|^{-1/3} \,\,\, (\gg 1).
\end{equation}
Following Titchmarsh, we divide the  double sum \eqref{exp} into $O(|t|^2\rho^2 X^{-4})$  subsums   where both variables run over intervals of length $\ll X^3 (|t| \| \bm \mu \|)^{-1}$. The crucial point is now that 
\begin{equation}\label{hess}
\| \text{Hess}(g(x_1, x_2))\|  \ll \frac{|t| \|\bm\mu \|}{X^3}, \quad \det \text{Hess}(g(x_1, x_2)) \gg \left(\frac{|t| \| \bm\mu \|}{X^3}\right)^2.
\end{equation}
The first bound is easy to see by direct calculation and trivial estimates as in \cite[p.\ 496]{Ti}, for the second bound we slightly deviate from Titchmarsh's argument on \cite[p.\ 497-498]{Ti} to avoid messy computations and use a Taylor approximation (since $\|{\bm  \mu} \| \leq \rho$ is much smaller than $\| \textbf{x} \| \asymp X$) getting
\begin{displaymath}
\begin{split}
\det \text{Hess}(g(x_1, x_2)) &=  \det \text{Hess} \left(-\mu_1f^{(1, 0)}(x_1, x_2) - \mu_2f^{(0, 1)}(x_1, x_2)  + h(x_1, x_2; \mu_1, \mu_2) \right),
\end{split}
\end{displaymath}
where $$\frac{\partial^i}{\partial x_1^i}\frac{\partial^j}{\partial x_2^j}  h(x_1, x_2; \mu_1, \mu_2) \ll   \frac{|t| \| \bm\mu \|^2}{\| \textbf{x} \|^{2+i+j}}$$
for $0 \leq i, j \leq 2$ (recall \eqref{norm}).  This gives 
\begin{displaymath}
\begin{split}
 \det \text{Hess}(g(x_1, x_2)) &=  \det \text{Hess} \Bigl(-\mu_1f^{(1, 0)}(x_1, x_2) - \mu_2f^{(0, 1)}(x_1, x_2) \Bigr)\left(1 + O\left(\frac{\|{\bm  \mu} \|}{ \| \textbf{x} \|}\right)\right) \\
&= \frac{4t^2\cdot \text{disc} (Q_{\pm}) \cdot Q_{\pm}({\bm\mu})}{Q_{\pm}(\textbf{x})^3} \left(1 + O\left(\frac{\|{\bm  \mu} \|}{ \| \textbf{x} \|}\right)\right),
\end{split}
\end{displaymath}
which implies in particular  the desired lower bound. Having \eqref{hess} available, a two-dimensional second derivative test \cite[Lemma $\delta$]{Ti}   bounds each of the subsums in \eqref{exp} by  $O(X^3 \| \bm \mu \|^{-1} |t|^{-1+\varepsilon})$ and hence  
$$\mathcal{S}_1 \ll \frac{\rho^2 |t|^{1+\varepsilon}}{\| \bm \mu \| X}.$$
 Now  Weyl differencing  \cite[Lemma $\beta$]{Ti}  shows 
$$\mathcal{S} \ll X^{-1}\Biggl(\frac{X^2}{\rho} + \frac{X}{\rho} \Bigl(\sum_{ 0 \not=  \| \bm \mu \| \ll \rho } \frac{\rho^2 |t|^{1+\varepsilon}}{\| \bm \mu \| X} \Bigr)^{1/2}\Biggr) \ll \frac{X}{\rho} +   \frac{\rho^{1/2} |t|^{1/2+\varepsilon}}{X^{1/2}} \ll |t|^{1/3+\varepsilon}$$
with our choice \eqref{rho} of $\rho$, which implies \eqref{toprove}. 

\section{Proof of Theorem \ref{thm2}}\label{sec3}

Let $k \geq 4$ be fixed. The theta-series
$$\theta_k(z) := \sum_{\textbf{x} \in \Bbb{Z}^k } e^{2\pi i z \| \textbf{x} \|_2^2} = \sum_{n=0}^{\infty} r_k(n) e(nz), \quad z \in \Bbb{H}_2,$$
with $r_k(n) = \#\{\textbf{x} \in \Bbb{Z}^k : \| x \|^2_2 = n\}$ 
is a modular form of weight $k/2$ for the group $\Gamma_1(4)$, and 
$$E(I_k, s) = \sum_{n=1}^{\infty} \frac{r_k(n)}{n^{ks/2}}$$
is the corresponding normalized $L$-series (the critical strip is $0 \leq \Re s \leq 1$). 
\subsection{The cuspidal contribution} 
Let $S(z) = \sum_n a(n) e(nz) \in S_{k/2}(\Gamma_1(4))$, say,  denote the orthogonal projection of $\theta$ onto the space of cusp forms. As is well-known \cite{Si}, this is 0 unless $k \geq 9$.  It can be decomposed into a sum of two cusp forms $S_+(z) + S_-(z)$, each of which is an eigenfunction of the level 4 Fricke involution.  The corresponding Dirichlet series $E_1(I_k, s) = \sum_n a(n) n^{-ks/2}$ is a linear combination of $L$-functions
$$\sum_{\pm}  L\left(\frac{ks}{2} - \frac{k/2 - 1}{2}, S_{\pm}\right).$$
These $L$-functions have no Euler product in general, but they are absolutely convergent in $\Re s > 1$ (by Hecke's mean square bound, since $S_{\pm}$ are cusp forms) and have holomorphic continuation together with  a functional equation (inherited from the Fricke involution). Hence both in the integral and half-integral weight case the convexity bound is available and gives 
\begin{equation}\label{cusp}
E_1(I_k, 1/2 + it) = \sum_{\pm}  L\left(1/2 + i kt/2, S_{\pm}\right) \ll (1 + |t|)^{1/2+\varepsilon}.
\end{equation}

We now consider the Eisenstein contribution $E_2(I_k, s) = E(I_k, s) - E_1(I_k, s)$ and distinguish several cases. 

\subsection{The case $k \equiv 0 \, (\text{mod } 4)$}
Here  the weight $k/2$ is even, and it follows from \cite[Theorem 44]{He} that  $$E_2(I_k, s) = c_k \zeta(ks/2) \zeta(ks/2-k/2+1) P(2^{-ks/2})$$
for some  polynomial $P$ with $P(0) = 1$ and some constant $c_k > 0$. The polynomial is given explicitly in \cite[3.8]{Sh},  
but we do not need this information. 
By the functional equation and Stirling's formula we have
$$  E_2\left(I_k, \frac{1}{2} + it\right) =  \zeta\left(\frac{k}{4} + \frac{ikt}{2}\right) \zeta\left(\frac{k}{4} - \frac{ikt}{2}\right)   P(2^{-k/4-itk/2})  \frac{\Gamma_{\Bbb{R}}( k/4 +itk/2 )}{\Gamma_{\Bbb{R}}( 1 - k/4 - itk/2  )}      \gg |t|^{k/4 - 1/2}$$
for $k > 4$ and $|t| \gg 1$. Together with \eqref{cusp} this gives the desired lower bound. The case $k=4$ is slightly more delicate. First  we recall that $E_1(I_4, 1/2 + it)  = 0$. Next, the zeta-function is now at the edge of the critical strip, and here we conclude from 
$$\int_{-T}^T |\zeta(1 + it)|^2  dt \sim 2 \zeta(2) T$$
(see \cite{BIR}) that $\zeta(1+it) = \Omega(1)$.  This gives again $E(I_4, 1/2 + it) = \Omega(|t|^{1/2})$, as desired. 

\subsection{The case $k \equiv 2 \, (\text{mod } 4)$}
Here the weight $k/2$ is odd. Again by \cite[Theorem 44]{He} we have
\begin{displaymath}
\begin{split}
E_2(I_k, s) = c_k&\left(\zeta(ks/2) L(ks/2 - k/2 + 1, \chi_{-4}) P_1(2^{-ks/2}) \right.\\
&\left.+ L(ks/2, \chi_{-4}) \zeta(ks/2 - k/2 + 1) P_2(2^{-ks/2})\right),
\end{split}
\end{displaymath}
and   it is convenient to record the exact description of $P_1$, $P_2$ from \cite[3.8]{Sh} as 
$$P_1(x) = \chi_{-4}(k/2), \quad P_2(x) = 2^{k/2 - 1} $$
(in particular $P_1$ and $P_2$ are constant). We apply  the functional equation for the $L$-series in the first term and for the zeta-function in the second term (noting that this involves slightly different $\Gamma$-factors) getting
$$ E_2(I_k, 1/2 + it) =  c_k \frac{ \Gamma_{\Bbb{R}}( k/4 - itk/2 ) }{\Gamma_{\Bbb{R}}( 1 - k/4 + itk/2 )} \left( A(it) + B(it)\right)$$
where
\begin{displaymath}
\begin{split}
 A(it) & = \zeta\left(\frac{k}{4} +\frac{ itk}{2}\right) L\left(\frac{k}{4} - \frac{itk}{2}, \chi_{-4}\right) P_1(2^{-k/4 - itk/2}) 4^{\frac{k}{4} - \frac{1}{2} - \frac{itk}{2}} \cot\left(\frac{\pi}{2}\left(1 - \frac{k}{4} + \frac{itk}{2}\right)\right)\\
 & = 2^{k/2-1}\zeta\left(\frac{k}{4} +\frac{ itk}{2}\right) L\left(\frac{k}{4} - \frac{itk}{2}, \chi_{-4}\right)\chi_{-4}(k/2)4^{-itk/2} + O(e^{-2|t|}), 
\\
 B(it) &  = 2^{k/2 - 1}L\left(\frac{k}{4} + \frac{itk}{2}, \chi_{-4}\right) \zeta\left(\frac{k}{4} - \frac{itk}{2}\right). \end{split}
\end{displaymath}
 The Gamma ratio is $\gg |t|^{k/4 - 1/2}$, and it remains to show that  $A(it) + B(it) = \Omega(1)$. To this end we consider
 \begin{equation}\label{consider}
 \int_{-T}^T 
 |A(it) + B(it)|^2 dt \gg T - 2\Re \int_{-T}^T  A(it)\overline{B(it)} dt. 
 \end{equation}
The integral on the right hand side equals
\begin{displaymath}
\begin{split}
&  2^{k-2} \chi_{-4}(k/2) \int_{-T}^T     \sum_{n_1, n_2, m_1, m_2} \frac{\chi_{-4}(n_2m_1)  }{(n_1n_2m_1m_2)^{k/4} }  \left(\frac{4n_1m_2}{n_2m_1}\right)^{-itk/2} dt + O(1)\\
 \ll\,& 1+\sum_{\substack{n_1, n_2, m_1, m_2\\ n_2m_1 \text{ odd}}} \frac{ 1  }{(n_1n_2m_1m_2)^{k/4} }  \min\left(T,  \Bigl|\log \frac{4n_1m_2}{n_2m_1}\Bigr|^{-1}\right).
  \end{split}
  \end{displaymath}
Since $4n_1m_2$ is even, but $n_2m_1$ is odd, we have
$$\min\left(T,  \Bigl|\log \frac{4n_1m_2}{n_2m_1}\Bigr|^{-1}\right) \leq T^{1/10}   \Bigl|\log \frac{4n_1m_2}{n_2m_1}\Bigr|^{-9/10}  \ll T^{1/10} \sqrt{n_1n_2m_1m_2}^{9/10}.$$ Since $k \geq 6$, the entire expression is $\ll T^{1/10}$, so that \eqref{consider} is $\gg T$, and hence $A(it) + B(it) = \Omega(1)$.

\subsection{The case $k$ odd}
Here it is most convenient to use the representation of $E_2(I_k, s)$ by Cohen \cite{Co}  $$E_2(I_k, s) = \sum_{D} \frac{L(\alpha, \chi_D) \zeta(ks) \zeta(ks + 2\alpha-1)}{D^s L(ks + \alpha, \chi_D)}, \quad \alpha = (3-k)/2,$$
where   the sum is over all fundamental discriminants $D$ with $(-1)^{(k-1)/2} D > 0$.  We apply the functional equation for the $L$-function in the numerator. This gives us 
$$E_2(I_k, s) = c_k \sum_{D} \frac{L(1-\alpha, \chi_D) \zeta(ks) \zeta(ks + 2\alpha-1)}{D^{ks/2+\alpha - 1/2} L(ks + \alpha, \chi_D)}$$
for some non-zero constant $c_k$. 
Since
$$\frac{\zeta(ks+2\alpha - 1)}{L(ks+\alpha, \chi_D)} = \sum_d \frac{1}{d^{ks+2\alpha-1}} \prod_{p \mid d} \left(1 - \frac{\chi_D(p)}{p^{1-\alpha}}\right),$$
we can rewrite this as
$$E_2(I_k, s) = c_k\zeta(ks )\sum_{\Delta} \frac{L(1-\alpha, \chi_{\Delta})}{\Delta^{s+\alpha - 1/2}}$$ 
where now the sum is over all (not necessarily fundamental)    discriminants  $\Delta$ with $(-1)^{(k-1)/2} \Delta > 0$. Let us write 
\begin{equation}\label{beta}
\beta = 1-\alpha = \frac{k-1}{2}, \quad w = \frac{ks}{2} + \alpha - 1/2 = \frac{ks}{2} - \frac{k}{2}  +1.
\end{equation}
With this notation, 
$$E_2(I_k, s) = c_k\zeta(2\beta + 2w - 1)\sum_{\Delta} \frac{L(\beta, \chi_{\Delta})}{\Delta^{w}}.$$
This can be expressed in terms of   the double Dirichlet series $$Z(s, w, \psi, \psi') = \zeta^{(2)}(2s +2w-1) \sum_{d \text{ odd}} \frac{L^{(2)}(s, \chi_d\psi)\psi'(d)}{d^w}$$ 
(where $\psi, \psi'$ are generally characters of conductor dividing  8 and  the Euler factors at 2 are removed),  which was 
investigated in detail in \cite{Bl}. For notational simplicity let us assume that $k \equiv 1$ (mod 4), so that $\Delta > 0$; the case $k \equiv 3$ (mod 4) is analogous. Every discriminant $\Delta$ has   precisely one of the forms $d 4^j$ with  $d \equiv 1$ (mod 4), $j \geq 0$ or   $d \equiv 2$ (mod 4), $j \geq 1$ or  $d \equiv 3$ (mod 4), $j \geq 1$.   Hence
\begin{displaymath}
\begin{split}
\sum_{\Delta} \frac{L(\beta, \chi_{\Delta})}{\Delta^{w}} & =  \sum_{d \equiv 1 \, (\text{mod }4)}  \frac{L(\beta, \chi_{d})}{d^w} + \sum_{d \equiv 1 \, (\text{mod }4)}  \frac{L^{(2)}(\beta, \chi_{d})}{d^w}  \frac{1}{4^w - 1} \\
& +\sum_{d \equiv 1 \, (\text{mod }2)}  \frac{L(\beta, \chi_{8d})}{(8d)^w}  \frac{4^{w}}{4^w - 1} +  \sum_{d \equiv 3 \, (\text{mod }4)}  \frac{L(\beta, \chi_{d})}{d^w}  \frac{1}{4^w - 1}, 
\end{split}
\end{displaymath}
so that
\begin{displaymath}
\begin{split}
  E_2(I_k, s) =  \frac{c_k}{2^w(4^w-1)(2^{2w + 2\beta - 1} - 1)} \sum_{\psi, \psi'} P_{\psi, \psi'}(2^w) Z(\beta, w, \psi, \psi') 
\end{split}
\end{displaymath}
for certain polynomials $P_{\psi, \psi'}$ (depending on $\beta$). 
The 16-dimensional vector of double Dirichlet series $Z(\beta, w, \psi, \psi')$ with $\psi, \psi'$ ranging over characters of conductor dividing 8 has a group of functional equations generated by $A : (\beta, w) \rightarrow (w, \beta)$ and $B : (\beta, w) \rightarrow (1 - \beta, s + \beta - 1)$. The corresponding 16-by-16  scattering matrices are given explicitly in \cite[(32), (33)]{Bl}. 
Applying the functional equation $B \circ A$, we obtain
$$E_2(I_k, s) =  C(\beta, w) \frac{\Gamma_{\Bbb{R}}(1-w)}{\Gamma_{\Bbb{R}}(w)}\sum_{\psi, \psi'} Z(1-w, \beta + w - 1/2, \psi, \psi') \tilde{P}_{\psi, \psi'}(2^w) + O(e^{-|\Im w|})$$
for certain polynomials $\tilde{P}_{\psi, \psi'}$ (not all 0) and a   factor $C(\beta, w)$ that on fixed vertical lines $\Re \beta, \Re w > 1$ is bounded from above and below.  As in the previous subsection, the error term comes from approximating the Gamma-factor $ \Gamma_{\Bbb{R}}(2-w)/\Gamma_{\Bbb{R}}(w+1)$ for odd characters  by $ \Gamma_{\Bbb{R}}(1-w)/\Gamma_{\Bbb{R}}(w)$. At $s = 1/2 + it$ with $|t| \geq 1$ and with $\beta$ and $w$ as in \eqref{beta}, we obtain
$$E_2(I_k, 1/2 + it) \gg |t|^{k/4 - 1/2} \Biggl| \sum_{\psi, \psi'} Z\left(\frac{k}{4} - \frac{itk}{2}, \frac{k}{4} + \frac{itk}{2}, \psi, \psi'\right) Q_{\psi, \psi'}(2^{itk/2})\Biggr| + O(e^{-|t|}) $$
for certain polynomials $Q_{\psi, \psi'}$ (not all 0), say
$$Q_{\psi, \psi'}(x) = \sum_{j=0}^J a_{\psi, \psi', j}x^j.$$ It remains to show that
\begin{equation}\label{remains}
\int_{-T}^T \Biggl| \sum_{\psi, \psi'} Z\left(\frac{k}{4} - \frac{itk}{2}, \frac{k}{4} + \frac{itk}{2}, \psi, \psi'\right) Q_{\psi, \psi'}(2^{itk/2})\Biggr|^2 dt  \gg T.
\end{equation}
The left hand side equals $\zeta^{(2)}(k - 1)^2$ times 
\begin{equation*}
\begin{split}
& \int_{-T}^T\Biggl| \sum_{n, d \text{ odd}} \sum_{j=1}^J\frac{\chi_d(n)}{(dn)^{k/4} }  \left(\frac{2^jn}{d}\right)^{ikt/2} \sum_{\psi, \psi'} \psi(n) \psi'(d) a_{\psi, \psi', j}    \Biggr|^2 dt \\
&= \sum_{n, n', d, d' \text{ odd}} \sum_{j, j' = 1}^J \int_{-T}^T  \left(\frac{2^j n d'}{2^{j'}n'd} \right)^{ikt/2} dt \sum_{\psi, \psi', \chi, \chi'}  \frac{\psi(n) \psi'(d) \chi(n') \chi'(d') \chi_d(n)\chi_{d'}(n') a_{\psi, \psi', j} \overline{a_{\chi, \chi', j'}}}{(dnd'n')^{k/4}}. 
\end{split}
\end{equation*} 
For the off-diagonal contribution $2^j n d' \not= 2^{j'}n' d$, we estimate the $t$-integral by
 $$  \min\left(T, \Bigl|\log \frac{2^j n d'}{2^{j'}n'd} \Bigr|^{-1}\right) \ll T^{3/5}  \Bigl|\log \frac{2^j n d'}{2^{j'}n'd} \Bigr|^{-2/5} \ll T^{3/5} \sqrt{nd'n'd'}^{2/5}.$$
 Since $k \geq 5/4$, the off-diagonal contribution is absolutely convergent and bounded by $T^{3/5}$. 
Since the character forces $(n, d) = (n', d') = 1$, the diagonal contribution is $d= d'$, $n = n'$, $j = j'$, and so  the left hand side of \eqref{remains} is $\zeta^{(2)}(k - 1)^2$ times 
$$ 2T  \sum_{\substack{n, d \text{ odd}\\ (n, d) = 1}}\sum_{j=1}^J \frac{1}{(nd)^{k/2}}  \Bigl| \sum_{\psi, \psi'} \psi(n) \psi'(d) a_{\psi, \psi', j} \Bigr|^2+ O(T^{3/5}). $$
Viewing $a_{\psi, \psi', j}$ for fixed $j$ as a function on the group $(\Bbb{Z}/2\Bbb{Z})^2 \times (\Bbb{Z}/2\Bbb{Z})^2$, it is clear that the $\psi, \psi'$-sum  can only vanish for all primitive residue classes $n, d$ modulo 8 if all $a_{\psi, \psi', j}$ vanish, which we have excluded. This proves \eqref{remains} and completes the proof of Theorem \ref{thm2}. 
  
\section{Proof of Proposition \ref{thm3}}\label{3} 

Let $D$ be a fixed   compact domain   of positive diagonal   matrices $\Xi = \text{diag}(\xi_1, \ldots, \xi_k)$. The key input is the following result of Jarn\'ik in \cite[Hilfssatz 6]{Ja} with
$$\sigma = k, \quad r_1 = \ldots = r_{\sigma} = 4, \quad x = T.$$
\begin{lemma}\label{Jarnik} There is a subset $M \subseteq D$ of full measure such that  for any $\Xi \in M$ and its associated quadratic form $Q(\textbf{x}) = \textbf{x}^{\top} \Xi \textbf{x}$ the following holds: let $T \geq 1$, $ T^{-1/2} \ll R \ll T^{4k}$.   Then
$$\int_R^{2R} |\theta(Q_{\pm}, 1/T + i\tau)|^4 d\tau  \ll T^{k+\varepsilon} R.$$ 
\end{lemma}
Here we used the same notation as before: $Q_+ = Q$ and $Q_-(\textbf{x}) = \textbf{x}^{\top}\Xi^{-1} \textbf{x}$. \\

We are now prepared to prove Proposition \ref{thm3}. We start from the approximate functional equation \eqref{afe}. Let $W$ be a fixed smooth function with support in $[1, 2]$. For the proof of \eqref{meansquare} it suffices (by the Cauchy-Schwarz inequality)  to show that
\begin{equation}\label{toshow}
 T^{-k/2} \int_{\Bbb{R}} W\left(\frac{t}{X}\right)  \Biggl|\sum_{\textbf{x} \not= 0} \frac{V_{T}(Q_{\pm}(\textbf{x}))}{Q_{\pm}(\textbf{x})^{\pm ikt/2 + iv}} \Biggr|^2  dt \ll X^{1+\varepsilon}
 \end{equation}
uniformly in $v \ll X^{\varepsilon}$, $T \ll X^{1+\varepsilon}$. Inserting \eqref{laplace}, the left hand side equals
\begin{equation}\label{equals}
T^{-k/2} \int_{(1/T)} \int_{(1/T)} \theta(Q_{\pm}, s_1) \overline{\theta(Q_{\pm}, s_2)} \int_{\Bbb{R}} W\left(\frac{t}{X}\right)  \mathcal{V}_{T, v}(s_1; t) \overline{\mathcal{V}_{T, v}(s_2; t) } dt\, \frac{ds_1\, ds_2}{(2\pi i)^2}.
\end{equation}
 Again we could compute $\mathcal{V}_{T, v}(s_1; t)$ by stationary phase uniformly in $t$ and then integrate over $t$, but an elementary argument suffices. As in Section \ref{sec2} we write $s_j = 1/T + i\tau_j$,  and we recall that $\mathcal{V}_{T, v}(1/T + i\tau; t) \ll |t|^{-A} \ll X^{-A}$  unless $|\tau| \asymp X/T$. 
 The $t$-integral equals
 \begin{displaymath}
 \begin{split}
 & \int_{\Bbb{R}^2} V_T(x) \overline{V_T(y)} \left(\frac{x}{y}\right)^{-iv} e^{s_1x + \overline{s_2}y}  \int_{\Bbb{R}}W\left(\frac{t}{X}\right)  \left(\frac{x}{y}\right)^{\mp ikt/2} dt  \, dx\, dy,\\
 & \int_{\Bbb{R}^2} V_T(xy) \overline{V_T(y)} x^{-iv} e^{\frac{1}{T}(x+1)y} e^{i y (\tau_1 x - \tau_2) }  \int_{\Bbb{R}}W\left(\frac{t}{X}\right)  x^{\mp ikt/2} dt  \, dx\, y\,dy.
 \end{split}
 \end{displaymath}
 Integrating by parts, we see that the $t$-integral is $\ll X^{-A}$ unless $x = 1 + O(X^{\varepsilon - 1})$, and in this range the  $y$-integral is  $\ll T^{-A}$ unless $\tau_1 x - \tau_2 \ll T^{\varepsilon-1}$, so that $$\tau_1 - \tau_2 \ll T^{\varepsilon - 1} + \frac{X}{T} X^{\varepsilon - 1} \ll X^{\varepsilon} T^{-1}.$$
 We conclude that \eqref{equals} is  
\begin{displaymath}
\begin{split}
& \ll T^{-A} + T^{-k/2} \int_{\substack{\tau_1 - \tau_2 \ll X^{\varepsilon} T^{-1}\\ \tau_1, \tau_2 \asymp X/T} } |\theta(Q_{\pm}, 1/T + i\tau_1) \overline{\theta(Q_{\pm}, 1/T + i\tau_2)}|   T^2 X^{\varepsilon}  d\tau_1\, d\tau_2 \\
& \ll T^{-A} + T^{-k/2} X^{\varepsilon} T \int_{\tau \asymp X/T} |\theta(Q_{\pm}, 1/T + i\tau)|^2 d\tau
\end{split}
\end{displaymath}
(using the Cauchy-Schwarz inequality in the second step). If $X/T \ll T^{4k}$, the bound  \eqref{toshow} follows now from another application of the Cauchy-Schwarz inequality   and Lemma \ref{Jarnik}, otherwise we estimate trivially $\theta(Q_{\pm}, 1/T + it) \ll T^{k/2}$, getting again the upper bound $T^{k/2 + 1} X^{\varepsilon} \ll X$.   This completes the proof of \eqref{meansquare}. 
 
 The bound \eqref{indiv} follows either from \eqref{meansquare} and the functional equation  (for instance as in \cite[p.\ 63]{Goo}) or directly from \eqref{afe}, \eqref{laplace}, \eqref{v1} and Lemma \ref{Jarnik}.

\end{document}